\documentclass[letterpaper, 10pt, conference]{ieeeconf}      

\IEEEoverridecommandlockouts                              
\newcommand{\ud}{\,\mathrm{d}}

\usepackage{graphics} 
\usepackage{epsfig} 
\usepackage{mathptmx} 
\usepackage{times} 
\usepackage{amsmath} 
\usepackage{amssymb}  
\usepackage{mathtools}
\usepackage{tabularx}
\usepackage{multirow}
\usepackage{subfigure}  

\usepackage{ifpdf}
\ifpdf
    \usepackage{graphicx}
    \usepackage{epstopdf}
\else
    \usepackage{graphicx}
\fi

\newtheorem{definition}{Definition}
\newtheorem{theorem}{Theorem}

\newtheorem{remark}{Remark}
\newtheorem{proposition}{Proposition}

\usepackage{algorithm}
\usepackage{algorithmic}

\newlength{\noteWidth}
\long\def\notes#1{\ifinner
             {\tiny #1}
             \else
              \marginpar{\parbox[t]{\noteWidth}{\raggedright\tiny #1}}
               \fi}

\def\notes#1{\typeout{#1 !!!}}  

\makeatletter
\newcommand{\rom}[1]{\romannumeral #1}
\newcommand*{\Rom}[1]{\expandafter\@slowromancap\romannumeral #1@}
\makeatother

\newcounter{rmnum}

\newcounter{anum}

\def\Sec#1{Sec.~\ref{#1}}

\def\IEEEQEDclosed{\mbox{\rule[0pt]{1.3ex}{1.3ex}}}

\def\qed{\nobreak\hfill\IEEEQEDclosed}
\def\UZ{\underline{\mathcal{Z}}}

\DeclareMathOperator{\Tr}{Tr}

\def\Re{\mathbb{R}}

\def\dt{\ud t}

\def\Expect{{\sf E}}
\def\K{{\sf K}}

\def\clG{{\cal G}}

\def\clL{{\cal L}}

\def\clV{{\cal V}}
\def\clW{{\cal W}}

\def\Cinf{{C^{\infty}}}

\newcommand{\half}[1]{\frac{#1}{2}}
\newcommand{\lr}[2]{\langle #1, #2 \rangle}

\def\x{{x}}

\def\k{{\sf k}}

\def\UZ{\mathcal{Z}_t}

\def\grad{\text{grad}}

\def\det{\text{det}}

\title{\LARGE \bf
Feedback Particle Filter on Matrix Lie Groups
}

\author{Chi Zhang, Amirhossein Taghvaei and Prashant G. Mehta 
\thanks{Financial support from the NSF CMMI grants 1334987 and 1462773 is gratefully acknowledged.}
\thanks{C.~Zhang, A.~Taghvaei and P. G.~Mehta are with the Coordinated
  Science Laboratory and the Department of Mechanical Science and
  Engineering at the University of Illinois at Urbana-Champaign (UIUC)
{\tt\scriptsize \{czhang54; taghvae2; mehtapg\}@illinois.edu}}
 }

\begin{document}

\maketitle

\begin{abstract}

This paper is concerned with the problem of continuous-time nonlinear 
filtering for stochastic processes on a compact and connected matrix Lie group without boundary, 
e.g. $SO(n)$ and $SE(n)$, 
in the presence of real-valued observations. This problem is important to numerous applications
in attitude estimation, visual tracking and robotic localization. 
The main contribution of this paper is to derive the feedback particle filter (FPF) 
algorithm for this problem. In its general form, the FPF provides a coordinate-free 
description of the filter that furthermore satisfies the geometric constraints of the manifold. 
The particle dynamics are encapsulated in a Stratonovich stochastic differential 
equation that preserves the feedback structure of the original Euclidean FPF. 
Specific examples for $SO(2)$ and $SO(3)$ are provided to help illustrate 
the filter using the phase and the quaternion coordinates, respectively.

\end{abstract}

\section{Introduction}
\label{sec:intro}

\medskip

There has been an increasing interest in the nonlinear filtering community to explore geometric approaches 
for handling constrained systems. In many cases, the constraints are described by 
smooth Riemannian manifolds, in particular the Lie groups. Engineering applications 
of filtering on Lie groups include: 
(\rom{1}) attitude estimation of satellites or aircrafts \cite{bonnabel2007cdc, barrau2015TAC}; 
(\rom{2}) visual tracking of humans or objects \cite{choi2011robust, kwon2007particle}; 
and (\rom{3}) localization of mobile robots \cite{barrau2014cdc, wang2008ijrr}. 
In these applications, the Lie groups of interest are primarily the matrix groups such as 
the special orthogonal group $SO(2)$ or $SO(3)$ and the special Euclidean group $SE(3)$.

\smallskip

This paper considers the continuous-time nonlinear filtering problem for matrix Lie groups 
in the presence of real-valued observations. The objective is to
obtain a generalization of the feedback particle filter (FPF) (see \cite{Tao_TAC}) in this non-Euclidean setting.
FPF is a continuous-time filtering algorithm that 
extends the feedback structure of the Kalman filter to general nonlinear non-Gaussian
filtering problems. For application problems in the Euclidean space, 
evaluation and comparison of FPF against the conventional particle filter appears in 
\cite{berntorp2015, stano2014, Adam_IF13}.

\medskip

The contributions of this paper are as follows:

\smallskip

\noindent {$\bullet$ \bf Feedback particle filter for Lie groups.} 
The extension of the FPF for matrix Lie groups is derived. The particle dynamics, expressed 
in their Stratonovich form, respect the manifold constraints. 
Even in the manifold setting, the FPF is \rom{1}) shown to admit an error correction feedback structure,  
and \rom{2}) proved to be an exact algorithm. Exactness means that, in the limit of large number 
of particles, the empirical distribution of the particles exactly matches the 
posterior distribution.

\smallskip

\noindent {$\bullet$ \bf Poisson equation on Lie groups.} 
Numerical implementation of the FPF requires approximation of the solution 
of a linear Poisson equation. The equation is described for the Lie group in 
an intrinsic coordinate-free manner. For computational purposes, a Galerkin scheme is proposed 
to approximate the solution. 

\smallskip

\noindent {$\bullet$ \bf Algorithms.} 
Specific examples for $SO(2)$ and $SO(3)$ are worked out, including expressions for the filter 
and the Poisson equation using a canonical choice 
of coordinates\textemdash the phase coordinate for $SO(2)$ and the quaternion coordinate for $SO(3)$. 

\medskip

Filtering of stochastic processes in non-Euclidean spaces has a rich history; c.f., 
\cite{caines1985IMA}, \cite{duncan1977}. 
In recent years, the focus has been on computational approaches 
to approximate the solution. 
Such approaches have been developed, e.g., by extending the 
classical extended Kalman filter (EKF) to Riemannian manifolds. 
EKF-based extensions have appeared for both discrete-time \cite{barrau2015TAC, barrau2013cdc} 
and for continuous time settings \cite{bonnabel2007cdc, bourmaud2015EKF}. 
Deterministic nonlinear observers have also 
been considered for $SO(3)$ \cite{mahony2008TAC, wu2015hybrid, batista2014attitude}, 
$SE(3)$ \cite{hua2011observer}, as well as for systems with other types of symmetry and 
invariance properties
\cite{lageman2010gradient, bonnabel2009observer}. A closely related theme is the use of 
non-commutative harmonic analysis for characterizing error propagation for 
rigid bodies \cite{park2008kinematic, kim2015bayesian, lee2008global}. 
These algorithms have also been applied extensively, e.g., for attitude estimation 
\cite{de2014, zanetti2009norm, markley2003attitude}.
Non-parametric approaches such as the particle filter (PF) have also been developed for 
Riemannian manifolds \cite{chiuso2000}, with extensive applications to visual tracking and 
localization \cite{kwon2014pf, kwon2007particle, choi2011robust}. Typically, PF algorithms adopt
discrete-time description of the dynamics and are based on importance sampling. 
Closely related approaches, such as the
Rao-Blackwellized particle filter \cite{kwon2010RBPF, barrau2014cdc} and the 
unscented Kalman filter (UKF) \cite{hauberg2013UKF}, have also been investigated for the Lie groups.

\smallskip

The remainder of this paper is organized as follows: After a brief overview in 
\Sec{sec:LG_prelim} of relevant preliminaries for matrix Lie groups, 
the filtering problem is formulated in \Sec{sec:filtering_LG}. In \Sec{sec:FPF_LG}, 
the generalization of the FPF algorithm to matrix Lie groups is presented, including both theory and algorithms. 
In \Sec{sec:examples}, examples for $SO(2)$ and $SO(3)$ are discussed. 
All the proofs appear in the appendix. 

\newpage

\section{Preliminaries}
\label{sec:LG_prelim}

This section includes a brief review of matrix Lie groups. 
The intent is to fix the notation used in subsequent sections.

The {\em general linear group}, denoted as $GL(n;\Re)$, is the group of $n\times n$ invertible matrices, 
where the group operation is matrix multiplication. The identity element is the identity matrix, 
denoted as $I$.
A {\em matrix Lie group}, denoted as $G$, is a closed subgroup of $GL(n;\Re)$. 
The {\em Lie algebra} of $G$, denoted as $\clG$, is the set of matrices $V$
such that the matrix exponential, $\exp (V)$, is in $G$. $\clG$ is a vector space whose 
dimension, denoted as $d$, equals the dimension of the group. 
$\clG$ is equipped with an inner product, denoted as $\lr{\cdot}{\cdot}_\clG$, 
and an orthonormal basis $\{E_1,...,E_d\}$ with $\lr{E_i}{E_j}_\clG = \delta_{ij}$. 
The space of smooth real-valued functions $f:G\rightarrow\Re$ is
denoted as $\Cinf(G)$.
\medskip

{\em Example:} The special orthogonal group $SO(3)$ is the group of $3\times3$ 
matrices $R$ such that $RR^T=I$ and $\det(R)=1$. The Lie algebra $so(3)$ is the 3-dimensional 
vector space of skew-symmetric matrices. An inner product is 
$\lr{\Omega_1}{\Omega_2}_\clG = (1/2)\Tr(\Omega_1^T\Omega_2)$, for $\Omega_1, \Omega_2 \in so(3)$, 
and an orthonormal basis $\{E_1,E_2,E_3\}$ of $so(3)$ is given by,
\[
E_1 = 
\begin{bmatrix}
 0 & 0 & 0 \\
 0 & 0 & -1 \\
 0 & 1 & 0
\end{bmatrix},
E_2 = 
\begin{bmatrix}
 0 & 0 & 1 \\
 0 & 0 & 0 \\
 -1 & 0 & 0
\end{bmatrix},
E_3 = 
\begin{bmatrix}
 0 & -1 & 0 \\
 1 & 0 & 0 \\
 0 & 0 & 0
\end{bmatrix}.
\]
These matrices have the physical interpretation of generating rotations about the 
three canonical axes (denoted as $e_1,e_2,e_3$) in $\Re^3$. $\det(\cdot)$ and $\Tr(\cdot)$ denote the 
determinant and trace of a matrix. \qed

\medskip

The Lie algebra can be identified with the tangent space at the
identity matrix $I$, and can furthermore be used to construct a basis
$\{E_1^x,...,E_d^x\}$ for the tangent space at $\x\in G$, where $E_n^x
= \x\,E_n$ for $n=1,...,d$.  Therefore, a smooth vector field, denoted as $\clV$, is expressed as,
  $$ \clV(\x) = v_1(\x) \, E_1^x + \cdots + v_d(\x) \, E_d^x, $$
with $v_n(x) \in \Cinf(G)$ for $n=1,...,d$. We write,
\begin{equation}
 \clV = x\,V,
 \label{eq:vector_field}
\end{equation}
where $V(x) := v_1(x)\,E_1 + \cdots + v_d(x)\,E_d$ is an element of the Lie algebra $\clG$ 
for each $\x\in G$. 

With a slight abuse of notation, the action of the vector 
field $\clV$ on $\psi\in\Cinf(G)$ is denoted as,
\begin{equation}
 V\cdot f(\x) := \frac{\ud}{\ud t}\Big|_{t=0} f \big(\x \, \exp(tV(\x))\big).
 \label{eq:action}
\end{equation}

The coordinates of $V$ are denoted as $v(x):=[v_1(\x),...,v_d(\x)]$. The inner product 
$\lr{\cdot}{\cdot}_\clG$ induces an inner product of two vector fields,
\begin{equation*}
 \lr{\clV}{\clW}(\x) := \lr{V}{W}_\clG(x) = \sum_{n=1}^d v_n(\x)w_n(\x).
\end{equation*}

It is understood that $V:G\rightarrow\clG$, where the identification with a vector field is 
through \eqref{eq:vector_field}, and the action is defined in \eqref{eq:action}. 
We will not use the $^\vee$ or $^\wedge$ notation to move between $V$ and its coordinates $v$, as is 
customary in the Lie groups \cite{chirikjian_book, bourmaud2015EKF}. 
This is because the $^\wedge$ notation 
is reserved for expectation, consistent with its use in stochastic processes.

As an example, we define next the notation for the vector field $\grad(\phi)$ 
for $\phi\in\Cinf(G)$,
\begin{equation}
 \grad(\phi)(\x) = \x \, \K(\x),
 \label{eq:gradient}
\end{equation}
where $\K(\x) = E_1\cdot\phi(\x) \, E_1 + \cdots + E_d\cdot\phi(\x) \, E_d \in \clG$, and according to 
\eqref{eq:action}, 
$(E_n\cdot \phi)(\x) := \frac{\ud}{\ud t} \big |_{t=0} \phi\big(x \,
\exp(tE_n)\big)$ for $n=1,...,d$.
The vector field acts on a function $f\in\Cinf(G)$ as,
\begin{equation}
 \K\cdot f(\x) = \sum_{n=1}^d E_n\cdot\phi(\x) \, E_n\cdot f(\x) = \lr{\grad(\phi)}{\grad(f)}(\x). 
 \label{eq:grad_inner_product}
\end{equation}

Apart from smooth functions, we will also need to consider other types of function spaces: 
For a probability measure $\pi$ on $G$, $L^2(G;\pi)$ denotes the Hilbert space of 
functions on $G$ that satisfy $\pi(|f|^2) < \infty$; 
$H^1(G;\pi)$ denotes the Hilbert space of functions $f$ such that $f$ and 
$E_n\cdot f$ (defined in the weak sense) are all in $L^2(G;\pi)$.

\section{Nonlinear Filtering Problem on Lie Groups}
\label{sec:filtering_LG}

\subsection{Problem statement}

We consider the following continuous-time system 
evolving on a Lie group $G$ with real-valued observations:
\begin{subequations}
\begin{align}
 \ud X_t & = X_t \, V_0(X_t)\dt + X_t \, V_1\circ\ud B_t, \label{eq:dyn} \\
 \ud Z_t & = h(X_t)\dt + \ud W_{t}, \label{eq:obs}
\end{align}
\end{subequations}
where $X_t\in G$ is the state at time $t$, $Z_t\in\Re$ is the observation, 
$V_0:G\rightarrow\clG$, $V_1 \in \clG$, $\{{B_t}\}$ and $\{W_t\}$ are mutually 
independent real-valued standard Wiener processes, 
which are also independent of the initial state $X_0$. 
The $\circ$ before $\ud B_t$ indicates that the stochastic differential equation (SDE) 
\eqref{eq:dyn} is defined in its Stratonovich form.

\medskip

\begin{remark}
On a smooth manifold, SDEs are usually constructed 
in their Stratonovich form instead of their It\^{o} form. 
The former respects the intrinsic geometry of the manifold~\cite{oksendal2003},
while the latter requires special geometric structures~\cite{gliklikh2010}, and 
in general does not maintain the manifold constraint~\cite{chirikjian_book}.~\qed
\end{remark}

\medskip

The objective of the filtering problem is to compute the conditional distribution of $X_t$ given 
the history of observations $\UZ=\sigma(Z_s:s\leq t)$. 
The conditional distribution, denoted as $\pi_t^*$, acts on a function $f\in\Cinf(G)$ 
according to,
  $$ \pi_t^*(f) := \Expect[f(X_t) | \UZ]. $$
$\pi_t^*$ is referred to as the filtered estimate.

\subsection{Filtering equation}

The filtering equation describes the evolution of the conditional distribution $\pi_t^*$. 
For the system \eqref{eq:dyn} and \eqref{eq:obs}, 
the Kushner-Stratonovich (K-S) filtering equation is (see \cite{crisan2009book}),  
\begin{align}
 \pi_t^*(f) = & \pi_0^*(f) + \int_0^t \pi_s^*({\cal L}^{*} f)\ud s \, + \notag \\
 & \int_0^t \big(\pi_s^*(fh)-\pi_s^*(h)\pi_s^*(f) \big) \big(\ud Z_s - \pi_s^*(h)\ud s \big),
 \label{eq:K-S}
\end{align}
for any $f\in\Cinf(G)$, where the operator ${\cal L}^{*}$ is given by,
\begin{equation}
 {\cal L}^{*} f = V_0\cdot f + \frac{1}{2} V_1\cdot (V_1\cdot f).
 \label{eq:L_star}
\end{equation}

\section{Feedback Particle Filter on Lie Groups}
\label{sec:FPF_LG}

This section extends the FPF algorithm originally proposed in \cite{Tao_TAC} to matrix
Lie groups, with necessary modifications to the original framework to account 
for the manifold structure.

\subsection{Particle dynamics and control architecture}

The feedback particle filter on a matrix Lie group $G$ is a controlled system 
comprising of $N$ stochastic processes $\{X_t^i\}_{i=1}^N$ with $X_t^i \in G$. 
The particles are modeled by the Stratonovich SDE,
\begin{equation}
 \ud X_t^i = X_t^i \big(V_0(X_t^i)+ u(X_t^i,t)\big)\ud t + X_t^i \, V_1\circ\ud B_t^i
                    + X_t^i \, \K(X_t^i,t) \circ \ud Z_{t},
 \label{eq:particle_dyn_uK}
\end{equation}
where $u(x,t), ~\K(x,t) :G \times [0,\infty) \rightarrow \clG$ are 
called {\em control} and {\em gain} function, respectively. These functions need to be chosen. 
The coordinates of $u$ and $\K$ are denoted as $[u_1,...,u_d]$ and $[\k_1,...,\k_d]$, 
respectively. Admissibility requirement is imposed on $u$ and $\K$:

\medskip

\begin{definition}
 {\em (Admissible Input)}:
 The functions $u(x,t)$ and $\K(x,t)$ are {\em admissible} if, 
 for each $t\geq0$, they are $\UZ-$measurable and we have 
 $\Expect[\sum_n|u_n(X_t^i,t)|]<\infty$, $\Expect[\sum_n |\k_n(X_t^i,t)|^2]<\infty$. \qed
\end{definition}

\medskip

The conditional distribution of the particle $X_t^i$ given $\UZ$ is denoted by $\pi_t$, 
which acts on $f\in\Cinf(G)$ according to,
 $$ \pi_t(f) := \Expect[f(X_t^i) | \UZ]. $$
The evolution PDE for $\pi_t$ is given by the 
proposition below. The proof appears in Appendix
\ref{apdx:prop_evolution}.

\medskip

\begin{proposition}
 Consider the particles $X_t^i$ with dynamics described by \eqref{eq:particle_dyn_uK}. 
 The forward evolution equation of the conditional distribution $\pi_t$ is given by,  
 \begin{equation}
  \pi_t(f) = \pi_0(f) 
   + \int_0^t \pi_s( \mathcal{L}f) \ud s 
   + \int_0^t \pi_s( \K\cdot f) \ud Z_{s},
  \label{eq:evo_p}
 \end{equation}
for any $f\in\Cinf(G)$, where the operator ${\cal L}$ is,
\begin{equation}
 {\cal L}f = (V_0+u)\cdot f + \frac{1}{2} V_1\cdot(V_1\cdot f) + \frac{1}{2}\K\cdot (\K\cdot f). 
 \label{eq:L}
\end{equation}
\label{prop:evolution} 
\qed
\end{proposition}

\medskip

\noindent {\bf Problem statement:}
There are two types of conditional distributions:
\begin{itemize}
 \item $\pi_t^{*}$: The conditional dist. of $X_t$ given $\UZ$.
 \item $\pi_t$: The conditional dist. of $X_t^i$ given $\UZ$.
\end{itemize}
The functions $\{u(x,t), \K(x,t)\}$ are said to be {\em exact} if $\pi_t=\pi_t^*$ 
for all $t \geq0$. Thus, the objective is to choose $\{u,~\K\}$ such that, 
given $\pi_0=\pi_0^{*}$, the evolution of the two conditional 
distributions are identical (see \eqref{eq:K-S} and \eqref{eq:evo_p}).

\medskip

\noindent{\bf Solution:}
The FPF on Lie groups represents the following choice of the gain function $\K$ and 
the control function $u$:

\noindent {\em 1. Gain function:} The gain function is obtained by solving a Poisson equation. 
Specifically, at each time $t$, Let $\phi_t\in H^1(G;\pi)$ be the solution of:
\begin{equation}
 \pi_t\big( \lr{\grad(\phi_t)}{\grad(\psi)} \big) = \pi_t\big( (h-\hat{h}_t) \psi \big), 
 \label{eq:poisson}
\end{equation}
for all $\psi\in H^1(G;\pi)$, where $\hat{h}_t=\pi_t(h)$. 
The gain function $\K$ is then given by,
$\x \, \K(\x, t) = \grad(\phi_t)(\x)$. Noting that (see \eqref{eq:gradient}), 
  $$ \grad(\phi_t)(\x) = E_1\cdot\phi_t(\x) \, E_1^x + \cdots + E_d\cdot\phi_t(\x) \, E_d^x, $$
where recall $E_n^x = x \, E_n$, we have,
  $$ \K(\x,t) = \k_1(\x,t)E_1 + \cdots + \k_d(\x,t)E_d, $$
with coordinates,
\begin{equation}
 \k_n(\x,t) = E_n\cdot \phi_t(\x)~,~~\text{for}~ n=1,...,d.
 \label{eq:gain_phi}
\end{equation}

\medskip

\noindent {\em 2. Control function:} The function $u$ is obtained as,
\begin{equation}
  u(x,t) = -\frac{1}{2}\K(x,t)(h(x)+\hat{h}_t).
 \label{eq:optimal_u}
\end{equation}

\noindent{\bf Feedback particle filter:}
Using these choice of $u$ and $\K$, FPF has the following representation: 
\begin{align}
 \ud X_t^i = & ~X_t^i \, V_0(X_t^i)\ud t + X_t^i \, V_1 \circ \ud B_t^i ~+ \notag \\
  & ~X_t^i \, \K(X_t^i,t) \circ \Big(\ud Z-\frac{h(X_t^i)+\hat{h}_t}{2}\ud t\Big).
 \label{eq:particle_dyn}
\end{align}

The consistency between $\pi_t^{*}$ and $\pi_t$ is asserted in the following theorem. The proof 
is contained in appendix \ref{apdx:thm_consistency}.

\medskip

\begin{theorem}
 Let $\pi_t^{*}$ and $\pi_t$ satisfy the forward evolution equations 
 \eqref{eq:K-S} and \eqref{eq:evo_p}, respectively. Suppose that the gain function $\K(x,t)$ 
 obtained using \eqref{eq:gain_phi}, and the control function $u(x,t)$ obtained 
 using \eqref{eq:optimal_u} are admissible. Then, assume $\pi_0=\pi_0^{*}$, we have, 
 $$ \pi_t(f) = \pi_t^{*}(f),$$
 for all $t \geq0$ and all function $f\in\Cinf(G)$.
\label{thm:consistency} 
\qed
\end{theorem}

\medskip

\begin{remark}
The admissibility of the control input leans on the existence,
uniqueness and regularity of the solution $\phi_t$ of the Poisson
equation~\eqref{eq:poisson}.  For the Euclidean space, this theory is
developed in~\cite{Laugesen_2015,Tao_CDC12} based on spectral
estimates for the $\pi_t^*$.  Extensions of these estimates to the
manifold settings is a subject of the continuing work.~\qed
\end{remark}

\subsection{Galerkin approximation}

The Poisson equation \eqref{eq:poisson} needs to be solved at each time step.
A Galerkin method is presented below to obtain an approximate solution.
Since the time $t$ is fixed, the explicit dependence on $t$ is
suppressed in what follows.  So, $\pi_t$ is denoted as $\pi$; $\phi_t$
is denoted as $\phi$ etc.

The function $\phi(x)$ is approximated as,
 $$ \phi(x) = \sum_{l=1}^L \kappa_l \, \psi_l(x), $$
where $\{\psi_l\}_{l=1}^L$ are a given (assumed) set of {\em basis functions} on
the manifold $G$.  Using~\eqref{eq:gain_phi}, the coordinates of the gain function $\K$ are then given by,
\begin{equation}
 \k_n(\x) = \sum_{l=1}^L \kappa_l \, E_n\cdot \psi_l(\x).
 \label{eq:gain}
\end{equation}

The finite-dimensional approximation of the Poisson equation~\eqref{eq:poisson} 
is to choose coefficients $\{\kappa_l\}_{l=1}^L$ such that
\begin{equation}
  \sum_{l=1}^L \kappa_l \, \pi\big( \lr{\grad(\psi_l)}{\grad(\psi)} \big)
     = \pi \big( (h-\hat{h})\psi \big),
  \label{eq:FEM}
\end{equation}
for all $\psi \in \text{span}\{\psi_1,...,\psi_L\} \subset H^1(G;\pi)$. 
On taking $\psi=\psi_1,...,\psi_L$, \eqref{eq:FEM} 
is compactly written as a linear matrix equation,
\begin{equation}
 A\kappa = b,
 \label{eq:kappa}
\end{equation}
where $\kappa := [\kappa_1,\hdots,\kappa_L]$ is a $L\times 1$ column
vector that needs to be computed.  The $L\times L$ matrix $A$ and the
$L\times 1$ vector $b$ are defined and approximated as,
\begin{align}
 [A]_{lm} & = \pi\big(\lr{\grad(\psi_l)}{\grad(\psi_m)} \big) \notag \\
 & \approx \frac{1}{N} \sum_{i=1}^N \lr{\grad(\psi_l)(X_t^i)}{\grad(\psi_m)(X_t^i)}, \notag \\
 & = \frac{1}{N} \sum_{i=1}^N \sum_{n=1}^d (E_n\cdot\psi_l)(X_t^i) \, (E_n\cdot\psi_m)(X_t^i), \label{eq:A} \\
 b_l & = \pi \big((h-\hat{h})\psi_l \big) 
 \approx \frac{1}{N}\sum_{i=1}^N (h(X_t^i)-\hat{h})\psi_l(X_t^i), \label{eq:b}
\end{align}
where $\hat{h}\approx\frac{1}{N}\sum_{i=1}^N h(X_t^i)$.

\medskip

\begin{remark}
 The Galerkin method is completely adapted to the data.  That is, no explicit 
 computation of the distribution is ever required. 
 Instead, one only needs to evaluate a given set of basis functions at the 
 particles $X_t^i$. The choice of basis functions $\{\psi_l\}_{l=1}^L$
 depends upon the problem. The functions $E_n\cdot\psi_l$ can typically be computed in
 an offline fashion. This is illustrated with examples in the next section. ~\qed
\end{remark}

\section{Examples}
\label{sec:examples}

This section contains two examples to illustrate the construction and 
implementation of the feedback particle filter.

\subsection{FPF on $SO(2)$}

$SO(2)$ is a 1-dimensional Lie group of rotation matrices $R$ such that 
$RR^T = I$ and $\det(R)=1$. An arbitrary element is expressed as,
\[
 R = R(\theta) = 
 \begin{bmatrix}
  \cos(\theta) & -\sin(\theta) \\
  \sin(\theta) & \cos(\theta)
 \end{bmatrix},
\]
where $\theta\in S^1$ is defined as the {\em phase coordinate}. 

The general form of the nonlinear filtering problem on $SO(2)$ is: 
\begin{align}
 \ud R_t & = R_t \, \omega(R_t) E \ud t + R_t \, E\circ\ud B_t, \label{eq:dyn_SO2} \\
 \ud Z_t & = h(R_t)\ud t + \ud W_t, \label{eq:obs_SO2}
\end{align}
where $\omega(\cdot)$ and $h(\cdot)$ are given real-valued functions on $SO(2)$, 
$\{B_t\}, \{W_t\}$ are independent standard Wiener processes in $\Re$, and 
$E$ is a basis of the Lie algebra $so(2)$,
\[
 E = 
 \begin{bmatrix}
  0 & -1 \\
  1 & 0
 \end{bmatrix}.
\]
The function $\omega(\cdot)$ has physical interpretation of the (local) angular velocity. 

The gain function $\K(R,t)=\k(R,t)E$ is a matrix in $so(2)$ with the coordinate $\k(R,t)$, 
a real-valued function on $SO(2)$. By the identification of $SO(2)$ and $S^1$ in 
terms of the phase coordinate $\theta$, define $\k(\theta,t) = \k(R(\theta),t)$. 
Similarly, define $\omega(\theta)=\omega(R(\theta))$, $h(\theta)=h(R(\theta))$, and
$\phi(\theta)=\phi(R(\theta))$.

It is straightforward to see that the filter expressed in the phase coordinate is given by,
\begin{align}
 \ud \theta_t^i = & ~\omega(\theta_t^i) \ud t + \ud B_t^i \, + \notag \\
 & ~\k(\theta_t^i,t)\circ(\ud Z_t - \frac{h(\theta_t^i)+\hat{h}_t}{2}\ud t),
  ~\text{mod} ~2\pi. 
 \label{eq:particle_SO2_theta}
\end{align}
At each time $t$, $\k(\theta,t)$ is obtained by solving the boundary value problem with respect to the phase 
coordinate. We suppress dependence on $t$, and write $\k(\theta)$ for $\k(\theta,t)$. 
With a slight abuse of notation, 
the action of $E$ on a smooth function is (see \eqref{eq:action}), 
\begin{equation*}
 E\cdot\phi(\theta) = \frac{\partial \phi}{\partial \theta}(\theta),
\end{equation*}
and the Poisson equation \eqref{eq:poisson} is expressed as,
\begin{equation}
 \pi \big( (E\cdot\phi) \, (E\cdot\psi) \big) = \pi \big( (h-\hat{h}) \psi \big),
 \label{eq:poisson_SO2}
\end{equation}
and needs to hold for all $\psi \in H^1(S^1;\pi)$. If $\pi$ has a probability density function $p$ on $S^1$, 
then one can write \eqref{eq:poisson_SO2} as,
\begin{equation*}
 \int_{0}^{2\pi} \frac{\partial \phi}{\partial\theta}(\theta)
 \frac{\partial \psi}{\partial\theta}(\theta) p(\theta) \ud \theta 
 = \int_{0}^{2\pi} (h(\theta)-\hat{h}) \psi(\theta) p(\theta) \ud \theta.
\end{equation*}

The solution of the Poisson equation on $S^1$ is approximated using a Fourier series basis. 
In the simplest case, these are just the first Fourier modes, in which case,
\begin{equation}
 \phi(\theta) = \kappa_1 \sin(\theta) + \kappa_2 \cos(\theta),
 \label{eq:basis_SO2}
\end{equation}
leading to the following formula in light of \eqref{eq:kappa}:
  \[
  \frac{1}{N}
    \begin{bmatrix}
      \sum_{i}\cos^2(\theta_t^i) & -\sum_{i}\cos(\theta_t^i)\sin(\theta_t^i) \\
      -\sum_{i}\cos(\theta_t^i)\sin(\theta_t^i) & \sum_{i}\sin^2(\theta_t^i) 
    \end{bmatrix}
    \begin{bmatrix}
     \kappa_1 \\
     \kappa_2
    \end{bmatrix}
  \]
  \[
  = \frac{1}{N}
    \begin{bmatrix}
     \sum_{i}(h(\theta_t^i)-\hat{h})\sin(\theta_t^i) \\
     \sum_{i}(h(\theta_t^i)-\hat{h})\cos(\theta_t^i)
    \end{bmatrix}.
  \] 
Finally, the gain function is obtained as,
\begin{equation}
 \k(\theta) = \kappa_1 \cos(\theta) - \kappa_2 \sin(\theta).
 \label{eq:gain_SO2}
\end{equation}

The resulting algorithm appears in \cite{Adam_2013ACC} where it is referred to as 
a coupled oscillator FPF. The filter is applicable to the problem of gait 
estimation in locomotion systems \cite{Adam_2012ACC}.

\subsection{FPF on $SO(3)$}

$SO(3)$ is a 3-dimensional Lie group (see the example in \Sec{sec:LG_prelim} for notation).
The nonlinear filtering problem is,
\begin{align}
 \ud R_t & = R_t \, \Omega \ud t + R_t \, V_1\circ\ud B_t, \label{eq:dyn_SO3} \\
 \ud Z_t & = h(R_t)\ud t + \ud W_t, \label{eq:obs_SO3}
\end{align}
where $\Omega, ~V_1\in so(3)$, and we write $\Omega = \omega_1E_1+\omega_2E_2+\omega_3E_3$. 
The coordinates $\omega_i$ may in general depend on $R_t$.

The feedback particle filter is given by,
\begin{align}
 \ud R_t^i = & R_t^i \, \Omega \ud t + R_t^i \, V_1\circ\ud B_t + \notag \\
 & R_t^i \, \K(R_t^i,t)\circ\Big(\ud Z_t - \frac{h(R_t^i) + \hat{h}_t)}{2}\ud t\Big),
 \label{eq:particle_SO3}
\end{align}
where the gain function $\K(R,t)=\k_1(R,t)E_1+\k_2(R,t)E_2+\k_3(R,t)E_3$ 
is an element of the Lie algebra $so(3)$. The coordinates
$(\k_1,\k_2,\k_3)$ are obtained by solving a Poisson equation~\eqref{eq:poisson}. 
We propose the following as basis functions:
\begin{equation}
 \begin{aligned}
 \phi_1(R) & = \frac{1}{2}e_2^T (R-R^T) e_3, ~\phi_2(R) = \frac{1}{2}e_3^T (R-R^T) e_1, \\
 \phi_3(R) & = \frac{1}{2}e_1^T (R-R^T) e_2, ~\phi_4(R) = \frac{1}{2}(\Tr(R) - 1), 
 \label{eq:basis_SO3}
\end{aligned}
\end{equation}
where $e_1, e_2, e_3$ are the canonical basis of $\Re^3$. 
The action of the basis $E_1, E_2, E_3$ is easily computed and given in Table-\Rom{1}. 

\begin{table}[H]
 \label{table:lie_deriv}
 \caption{Action of $E_n$ on basis functions}
 \centering
 {\normalsize
 \begin{tabular}{c|c|c|c}
  \hline
             & $E_1\cdot$ & $E_2\cdot$ & $E_3\cdot$ \\ \hline
    $\phi_1$ & $-(R_{22}+R_{33})/2$ & $R_{21}/2$ & $R_{31}/2$ \\ \hline
    $\phi_2$ & $R_{12}/2$ & $-(R_{11}+R_{33})/2$ & $R_{32}/2$ \\ \hline
    $\phi_3$ & $R_{13}/2$ & $R_{23}/2$ & $-(R_{11}+R_{22})/2$ \\ \hline
    $\phi_4$ & $(R_{23}-R_{32})/2$ & $(R_{31}-R_{13})/2$ & $(R_{12}-R_{21})/2$ \\
  \hline
 \end{tabular}
 }
\end{table}

Using Table-\Rom{1}, the $4\times4$ matrix $A$ and the $4\times1$ vector $b$  
are assembled according to \eqref{eq:A} and \eqref{eq:b}, respectively.
The solution of the linear equation \eqref{eq:kappa} is a $4\times1$ vector, 
denoted as $\kappa=[\kappa_1,\kappa_2,\kappa_3,\kappa_4]$. 
Denoting  
\begin{equation}
  \varUpsilon = 
 \begin{bmatrix}
  \kappa_4 & -\kappa_3 & \kappa_2 \\
  \kappa_3 & \kappa_4 & -\kappa_1 \\
  -\kappa_2 & \kappa_1 & \kappa_4
 \end{bmatrix},
 \label{eq:varUpsilon}
\end{equation}
the coordinate functions of the gain have a succinct representation,
\begin{equation}
 \k_n(R) = \half{1}\Tr(RE_n\varUpsilon),~~\text{for}~n=1,2,3.
\end{equation}

For computational reasons, quaternions is a preferred choice for simulating rotations in $SO(3)$ 
\cite{markley2003attitude, jahanchahi2014class}. 
A unit quaternion has a general form,
  $$ q=\big( \cos(\frac{\theta}{2}), ~\sin(\frac{\theta}{2})\omega_1, ~\sin(\frac{\theta}{2})\omega_2, ~
  \sin(\frac{\theta}{2})\omega_3 \big)^T, $$
which represents rotation of angle $\theta$ about the axis defined 
by the unit vector $(\omega_1, \omega_2, \omega_3)^T$. A quaternion is also 
written as $q = (q_0,~q_1,~q_2,~q_3)^T$.

In the following, the FPF is described for the quaternion coordinates. In these coordinates, 
the four basis functions (counterparts of \eqref{eq:basis_SO3}) are,
\begin{equation}
 \begin{aligned}
 \phi_1(q) & = 2q_1q_0 , ~~\phi_2(q) = 2q_1q_0 , \\
 \phi_3(q) & = 2q_1q_0 , ~~\phi_4(q) = 2q_0^2-1 .  \\
 \end{aligned}
 \label{eq:basis_SO3_quat}
\end{equation}

In order to compute the matrix $A$ and the vector $b$, the formulae
for the action 
of $E_1,~E_2,~E_3$ on these basis functions appear in Table-\Rom{2}.

\begin{table}[H]
 \label{table:lie_deriv_quat}
 \caption{Action of $E_n$ on basis functions using quaternion}
 \centering
 {\normalsize
 \begin{tabular}{c|c|c|c}
  \hline
             & $E_1\cdot$       & $E_2\cdot$       & $E_3\cdot$      \\ \hline
    $\phi_1$ & $q_1^2-q_0^2$    & $q_1q_2-q_3q_0$  & $q_1q_3+q_2q_0$ \\ \hline
    $\phi_2$ & $q_1q_2+q_3q_0$  & $q_2^2-q_0^2$    & $q_2q_3-q_1q_0$ \\ \hline
    $\phi_3$ & $q_1q_3-q_2q_0$  & $q_2q_3+q_1q_0$  & $q_3^2-q_0^2$   \\ \hline
    $\phi_4$ & $2q_1q_0$        & $2q_2q_0$        & $2q_3q_0$       \\
  \hline
 \end{tabular}
 }
\end{table}

As before, the solution of the linear matrix equation is denoted as 
$\kappa=(\kappa_1,\kappa_2,\kappa_3,\kappa_4)$, and the coordinates of the 
gain function are obtained as, 
  $$ \k_n(q,t) = \frac{1}{2}\Tr(R(q)E_n\varUpsilon), $$
where $\varUpsilon$ is defined in \eqref{eq:varUpsilon}, and $R(q)$ is obtained 
using the conversion rule between rotation matrices and quaternions (see \cite{markley2003attitude}).
  
Finally, the filter in the quaternion coordinates has the following form,
\begin{equation}
 dq_t^i = \frac{1}{2}\Lambda\big(V(q_t^i)\big) \, q_t^i 
 + \frac{1}{2}\Lambda\big(\K(q_t^i,t)\big) \, q_t^i \circ (\ud Z_t-\frac{h(q_t^i)+\hat{h}_t}{2}\ud t),
 \label{eq:dyn_quat}
\end{equation}
where $\K(q_t^i,t),~V(q_t^i)\in so(3)$, $V(q_t^i) = \Omega\ud t + V_1\circ\ud B_t^i$
and the $4\times4$ matrix $\Lambda(\K)$ is given by,
\[
 \Lambda\big( \K \big) := 
 \begin{bmatrix}
  0 & -\k_1 & -\k_2 & -\k_3 \\
  \k_1 & 0 & \k_3 & -\k_2 \\
  \k_2 & -\k_3 & 0 & \k_1 \\
  \k_3 & \k_2 & -\k_1 & 0
 \end{bmatrix},
\]
and similarly for $\Lambda\big(V(q_t^i)\big)$.

\medskip

\begin{remark}
Consider the special case where the dynamics are restricted to the
subgroup $SO(2)$ of $SO(3)$. In this case, the filter \eqref{eq:dyn_quat} 
for $SO(3)$ reduces to the filter \eqref{eq:particle_SO2_theta} for $SO(2)$. 
To see this, note that with the axis of rotation $(\omega_1,~\omega_2,~\omega_3)$ fixed, 
the four basis functions are given by,
\begin{equation*}
 \begin{aligned}
 \phi_1(q) & = 2q_1q_0 = \sin(\theta)\omega_1, ~\phi_2(q) = 2q_1q_0 = \sin(\theta)\omega_2, \\
 \phi_3(q) & = 2q_1q_0 = \sin(\theta)\omega_3, ~\phi_4(q) = 2q_0^2-1  = \cos(\theta). \\
 \end{aligned}
 \label{eq:basis_SO3_quat_theta}
\end{equation*}
These functions span a 2-dimensional space, same as the Fourier basis functions $\{\sin(\theta), \cos(\theta)\}$ for the $SO(2)$ problem. \qed
\end{remark}

\section{Conclusions}
\label{sec:conclusion}

In this paper, the generalization of the feedback particle filter to the continuous-time 
filtering problem on matrix Lie groups was presented. The formulation was shown to respect the intrinsic 
geometry of the manifold and preserve the error correction-based feedback structure 
of the original FPF. Algorithms were described and illustrated with 
examples for $SO(2)$ and $SO(3)$.

The continuing research includes application and evaluation of the filter to attitude estimation and 
robot localization; comparison of the FPF with existing algorithms based on EKF 
and the particle filter; and extension of the FPF for filtering stochastic processes 
where the observation also evolves on manifold (c.f., \cite{caines1985IMA, said2013obs}).

\appendix
\label{sec:appendix}

\subsection{Proof of Proposition \ref{prop:evolution}}
\label{apdx:prop_evolution}

The solution $X_t^i$ to the 
Stratonovich SDE \eqref{eq:particle_dyn_uK} is a continuous semimartingale on the Lie group. 
Furthermore, for any smooth function $f: G \rightarrow \Re$, $f(X_t^i)$ is also a 
continuous semimartingale \cite{shigekawa1984}, satisfying, 
\begin{align}
 \ud f(X_t^i) = & (V_0+u)\cdot f(X_t^i)\ud t + V_1\cdot f(X_t^i) \circ \ud B_{t}^i \notag \\
 & + (\K\cdot f)(X_t^i) \circ \ud Z_{t}.
 \label{eq:df_strat}
\end{align}

To avoid the technical difficulty in taking expectation of the Stratonovich stochastic 
integrals, we convert \eqref{eq:df_strat} to its It\^{o} form using 
the formula given in \cite{watanabe1981}:
For continuous semi-martingales $X, Y, Z$,
\begin{align}
 & Y\circ \ud X = Y \ud X + \frac{1}{2}\ud X \ud Y, \label{eq:eq1} \\
 & (X\circ\ud Y) \ud Z = X (\ud Y \ud Z). \label{eq:eq2}
\end{align}

To convert the second term of the right hand side of \eqref{eq:df_strat}, take $Y$ in \eqref{eq:eq1} to be 
$(V_1\cdot f)(X_t^i)$ and $X$ to be $B_{t}^i$, we have,
\begin{equation}
 (V_1\cdot f)(X_t^i) \circ \ud B_{t}^i = (V_1\cdot f)(X_t^i) \ud B_{t}^i
 + \frac{1}{2} \ud (V_1\cdot f)(X_t^i) \ud B_{t}^i.
 \label{eq:convert_E}
\end{equation}
Then replace $f$ by $V_1\cdot f$ in \eqref{eq:df_strat} to obtain,
\begin{align*}
 \ud (V_1\cdot f) = & (V_0+u)\cdot (V_1\cdot f)\ud t + 
 V_1\cdot(V_1\cdot f) \circ \ud B_{t}^i \notag \\
 & + \K\cdot (V_1\cdot f) \circ \ud Z_{t}.
 \label{eq:dEf}
\end{align*}
Using \eqref{eq:eq2} and It\^{o}'s rule, 
\begin{equation}
 \ud (V_1\cdot f)(X_t^i) \ud B_{t}^i = V_1\cdot(V_1\cdot f)(X_t^i) \ud t,
\end{equation}
which when substituted in \eqref{eq:convert_E} yields,
 $$ V_1\cdot f(X_t^i) \circ \ud B_{t}^i = V_1 \cdot f(X_t^i) \ud B_{t}^i 
  + \frac{1}{2} V_1\cdot(V_1\cdot f)(X_t^i) \ud t. $$

The third term on the right hand side of \eqref{eq:df_strat} is similarly converted. 
The It\^{o} form of \eqref{eq:df_strat} is then given by,
\begin{equation*}
 \ud f(X_t^i) = \mathcal{L}f(X_t^i) \ud t + V_1 \cdot f(X_t^i) \ud B_{t}^i
 + (\K\cdot f)(X_t^i,t) \ud Z_{t},
\end{equation*}
where the operator $\mathcal{L}$ is defined by,
\begin{equation}
 \mathcal{L}f := (V_0+u)\cdot f + \frac{1}{2} V_1\cdot(V_1\cdot f) + \frac{1}{2}\K\cdot (\K\cdot f).
\end{equation}

The solution of $f(X_t^i)$ is obtained as, 
\begin{align}
 f(X_t^i) = f(X_0^i) & + \int_0^t \mathcal{L}f(X_s^i)\ud s 
            + \int_0^t V_1\cdot f(X_s^i) \ud B_{s}^i \notag \\
            & + \int_0^t (\K\cdot f)(X_s^i) \ud Z_{s}. \notag
\end{align}
By taking conditional expectation on both sides and interchanging 
expectation and integration,
 $$ \pi_t (f) = \pi_0(f) 
   + \int_0^t \pi_s( \mathcal{L}f ) \ud s 
   + \int_0^t \pi_s( \K\cdot f ) \ud Z_{s},
 $$
which is the desired formula \eqref{eq:evo_p}.

\subsection{Proof of Theorem \ref{thm:consistency}}
\label{apdx:thm_consistency}

Using~\eqref{eq:K-S} and~\eqref{eq:evo_p} and the expressions for the
operators $\clL^*$ and $\clL$, it suffices to show that
 \begin{align}
 \pi_s(u\cdot f) \ud s & + \frac{1}{2} \pi_s \big( \K\cdot(\K\cdot f) \big) + \pi_s(\K\cdot f) \ud Z_s \notag \\
 & = \big(\pi_s(fh)-\pi_s(h)\pi_s(f)\big) \big(\ud Z_s - \pi_s(h)\ud s \big),
 \label{eq:WTS}
 \end{align}
for all $0 \leq s \leq t$, and all $f\in\Cinf(G)$.

On taking $\psi=f$ in \eqref{eq:poisson} and using the
formula~\eqref{eq:grad_inner_product} for the inner product, 
\begin{equation}
 \pi_s(\K\cdot f) = \pi_s \big( (h-\pi_s(h))f \big).
 \label{eq:Kf}
\end{equation}

Using the expression~\eqref{eq:optimal_u} for the control function and noting $\hat{h}_s=\pi_s(h)$,
  $$ u\cdot f = -\frac{1}{2}(h-\pi_s(h)) \K\cdot f - \pi_s(h) \K\cdot f. $$
Using \eqref{eq:Kf} repeatedly then leads to, 
\begin{equation}
 \begin{aligned}
 \pi_s(u\cdot f) & = -\frac{1}{2} \pi_s \big( (h-\pi_s(h)) \K\cdot f \big) - \pi_s(h) \pi_s(\K\cdot f) \\
 & = -\frac{1}{2} \pi_s \big( \K\cdot(\K\cdot f) \big) - \pi_s(h) \pi_s \big( (h-\pi_s(h))f \big).
 \label{eq:uf}
\end{aligned}
\end{equation}

The desired equality~\eqref{eq:WTS} is now verified by substituting in \eqref{eq:Kf} and \eqref{eq:uf}.

\bibliographystyle{plain}
\bibliography{16ACC_LGFPF}

\end{document}